\documentclass[twocolumn]{amsart}
\usepackage[mathscr]{euscript}
\usepackage{amssymb}
\newcommand\C{\mathbb C}

\newcommand\K{\mathcal K}

\newcommand\M{\mathcal M_{\check{X\,}\!}}
\renewcommand\O{\mathcal O}
\newcommand\PP{\mathbb P}

\newcommand\Q{\mathbb Q}
\newcommand\R{\mathbb R}

\newcommand\X{\check{X\,}\!}
\newcommand\ppi{\check{\pi\,}\!}
\newcommand\Z{\mathbb Z}

\newcommand{\rt}[1]{\stackrel{#1\,}{\rightarrow}}
\newcommand{\Rt}[1]{\stackrel{#1\,}{\longrightarrow}}

\newcommand\res{\arrowvert^{}_}

\newcommand\ip{\lrcorner}

\newcommand\take{\backslash}
\newcommand{\mat}[4]{\left(\begin{array}{cc} \!\!#1 & #2\!\! \\ \!\!#3 &
#4\!\!\end{array}\right)}

\newcommand\rk{\operatorname{rank}}

\newcommand\im{\operatorname{im}}
\newcommand\id{\operatorname{id}}
\newcommand\vol{\operatorname{vol}}

\newcommand\Hom{\operatorname{Hom}}
\newcommand\Ext{\operatorname{Ext}}

\newcommand\Proj{\operatorname{Proj}\,}

\renewcommand\Re{\operatorname{Re}}

\newcommand\beq[1]{\begin{equation}\label{#1}}
\newcommand\eeq{\end{equation}}
\newcommand\beqa{\begin{eqnarray*}}
\newcommand\eeqa{\end{eqnarray*}}

\title{The geometry of mirror symmetry}
\author{Richard P. Thomas, Imperial College, London, UK}
\date{}

\begin{document}
\maketitle

\section{Introduction}

Mirror symmetry was discovered in the late 1980s by physicists studying superconformal
field theories (SCFTs). One way to produce SCFTs is from \emph{closed string theory};
in the Riemannian (rather than Lor\-entzian) theory the string's worldline gives a
map of a Riemannian 2-manifold into the target with an action which is
conformally invariant, so the 2-manifold can be thought of as a \emph{Riemann
surface} with a complex structure. Making sense of the infinities
in the quantum theory (supersymmetry and anomaly cancellation) forces the
target to be 10 dimensional -- Minkowski space times by a 6-manifold $X$
-- and $X$ to be (to first order) Ricci flat and so to have holonomy in
$SU(3)$. That is $X$ is a \emph{Calabi-Yau 3-fold} $(X,\Omega,\omega)$.
So SCFTs come from
\emph{$\sigma$-models} (mapping Riemann surfaces into Calabi-Yau
3-folds) but, it turns out, in two different ways -- the A-model and the
B-model. Deformations of the SCFT and either $\sigma$-model are isomorphic,
so over an open set the two coincide. So it was natural to conjecture that
pretty much all of the relevant SCFTs came from geometry -- from an A or
B $\sigma$-model.
In particular, the A-model of a Calabi-Yau $X$ should therefore give the
same SCFT as the B-model on another Calabi-Yau $\X$. It turns out then that
the A-model on $\X$ should also be isomorphic to the B-model on $X$; thus
mirror symmetry should give an involution on the set of Calabi-Yau manifolds. (The full picture is
slightly more complicated, to do with large complex structure limits, described
below, multiple mirrors and flops.) By studying the
SCFTs, Greene and Plesser predicted the mirror of the simplest Calabi-Yau 3-fold,
the quintic in $\PP^4$, and mirror symmetry was born.

Topological observables, i.e. certain path integrals over the space of all
maps, can be calculated by the semi-classical approximation as integrals
over the space of classical minima -- (anti)holomorphic curves in the Calabi-Yau
(these minimize volume in a fixed homology class). From the zero homology
class we get the constant maps -- points in $X$ -- and so integrals over
$X$. In some cases, by Poincar\'e duality, these can be thought of as intersections
of cycles; we think of the string worldsheet lying at a point of intersection.
When the worldsheet has a nontrivial homology class it allows more general
`intersections' where the cycles need not intersect but are connected by
a holomorphic curve, giving a perturbation of the usual intersection product
on cohomology called \emph{quantum cohomology}. Namely, there is a contribution
$(a.\beta)(b.\beta)(c.\beta)e^{\int_\beta\omega}$ to the quantum triple
product $a.b.c$ of three four-cycles $a,b,c\in H^{1,1}\cong H^2 \cong H_4$
from each holomorphic curve $\beta$ (of genus 0, in the 0-loop approximation
to the physics) in $X$ of
area $\int_\beta\omega$ (where $\omega$ is the K\"ahler form). The A-model
correlation functions can be determined from this data; the B-model computation
involves no such quantum correction and can be computed purely in terms of
integrals over cycles (``periods") and their derivatives, discussed in
(\ref{Yuk}) below. So it is in some
sense easier and, in a historic tour-de-force, was calculated in \cite{COGP}
for the Greene-Plesser mirror of the quintic. Comparing with the A-model
computation on the quintic gave remarkable predictions about the number of
holomorphic rational curves on the quintic. These
were way beyond mathematical capabilities at the time, and sparked enormous
mathematical interest. The predictions (and more) have now been proved
to be true by Givental and Lian-Liu-Yau, while mirror symmetry has begun
to be understood geometrically. But in some sense the
mathematical \emph{reason} for the relationship between the Yukawa couplings
and the quantum cohomology of the mirror is still a little mysterious; it
is the hardest
part of mirror symmetry to see in the geometry, yet for the physics it was
the easiest and first prediction.

We survey, non-chronologically, some of the geometry of mirror symmetry
as it is now understood, mainly in dimension $n=3$. For the many
topics omitted the reader should consult the references.

\section{The geometric set-up} \label{1}

A \emph{Calabi-Yau 3-fold} $(X,\Omega,\omega)$ is a K\"ahler manifold $(X,\omega)$
with a holomorphic trivialization $\Omega$ of its \emph{canonical bundle}
$$
K_X=\Lambda^3_\C T^*X
$$
(i.e. a nowhere vanishing holomorphic volume form, locally $dz_1\wedge dz_2\wedge
dz_3$), and $b_1(X)=0$. It follows that the Hodge
numbers $h^{0,2},\,h^{0,1}$ vanish, and so $H^2(X,\C)=H^{1,1}$ and
$H^3(X,\R)\cong H^{2,1}\oplus H^{3,0}$. By Yau's theorem the K\"ahler metric can be
changed within its $H^2(X,\R)$ cohomology class to a unique \emph{Ricci-flat}
K\"ahler metric; equivalently $\Omega$ is parallel, so the induced metric on $K_X$ is \emph{flat}.
Roughly speaking, mirror symmetry swaps the symplectic or K\"ahler structure
$\omega$ on $X$ with the complex structure (encoded in $\Omega$,
up to scaling by $\C^*$)
on the (conjectural) mirror $\X$. K\"ahler deformations are unobstructed,
forming an open set $\K_X$ in $H^2(X,\R)$. Its closure $\overline\K_X$
is sometimes extended by adding the K\"ahler cones of all birational models of
$X$ to give Kawamata's \emph{moveable cone}. This is because work of Aspinwall,
Greene, Morrison and Witten suggested
that all birational models of $X$ are indistinguishable to string theory
and so are all mirrors of $\X$, corresponding to a different choice of
$(1,1)$-form $\omega$ which is a K\"ahler form on one model only. $\K_X$
is also complexified by including in the A-model data any ``B-field"
$B\in H^2(X,\R/\Z)$, and divided by holomorphic automorphisms of $X$, to give a moduli space of \emph{complex} dimension
$h^{1,1}(X)$. Deformations of complex structure are also unobstructed
by the non-trivial Bogomolov-Tian-Todorov theorem, so form
a smooth space with tangent space
$$
H^1(T\X)\mathop{\Rt{\ip\,\Omega}}_\simeq H^1(\Lambda^2T^*\X)=H^{2,1}(\X).
$$
(Given a deformation of complex structure, this isomorphism takes the
$H^{2,1}$-component of the derivative of the $(3,0)$-form $\Omega$.)
So for the moduli spaces to match up, we get the first and simplest prediction
of mirror symmetry, that
\beq{hodge}
h^{1,1}(X)=h^{2,1}(\X) \quad\text{and}\quad h^{2,1}(X)=h^{1,1}(\X).
\eeq
This is where mirror symmetry gets its name, the above relation making
the Hodge diamonds of $X$ and $\X$ mirror images of each other.

As the complexified K\"ahler cone is a \emph{tube domain} it has natural
partial complex compactifications (due to Looijenga, and suggested in the
context of mirror symmetry by Morrison \cite{M}).
The simplest case is where we ignore the moveable cone and automorphisms
and assume there is an integral basis
$e_1,\ldots e_n$ of both $\overline\K_X$ and
$H^2(X,\Z)/$torsion. The complexified K\"ahler moduli space is then
$$
\K^\C_X:=H^2(X,\R)/H^2(X,\Z)+i\K_X=\{B+i\omega\},
$$
with natural coordinates $x_i,\,y_i\ge0$ pulled back from the first and second
factors respectively, induced by the $e_i$. $x_i$ is multiply-valued with
integer periods, so
\beq{q}
z_i=\exp(2\pi i(x_i+iy_i))
\eeq
is a well defined holomorphic coordinate, giving an isomorphism to the
product of $n$ punctured unit disks in $\C$:
$$
\K_X^\C\cong(\Delta^*)^n=\{(z_i)\colon0<|z_i|\le1\}\subset(\C^*)^n.
$$
The compactification $\Delta^n$ comes from adding in the origins in the disks,
which we reach by going to infinity (in various directions) in $\K^\C_X$.
We call the point $(0,\ldots,0)\in\Delta^n$ the \emph{large K\"ahler limit point}
(LKLP) in this case.
Moving along the ray generated by $\sum k_ie_i\in\K_X,\ k_i\ge0,$
complexifies in the holomorphic structure (\ref{q}) to give the analytic curve
\beq{ni}
z_i^{k_j}=z_j^{k_i},\ \ \forall\,i,j,
\eeq
in $\K^\C_X$.
For $k_i\in\Q\ \forall i$ this extends to a complete curve in the
compactification. Without loss of generality we can assume the $k_i$ are
integers with no common factor, then the link of the curve winds around the
LKLP
$(0,\ldots,0)\in\Delta^n$ with winding number
$$
(k_1,\ldots,k_n)\in\pi_1(H^2(X,\R)/H^2(X,\Z)+i\K_X)
$$$$
=H^2(X,\Z)=\Z.{e_1}\oplus\ldots\oplus\Z.{e_n}.
$$
This is because multiplying the ray $\R\,.\!\sum k_ie_i\in\K_X$ by $i$
gives the direction $\R\,.\!\sum k_ie_i$ in the space $H^2(X,\R)/H^2(X,\Z)$
of B-fields, with the given winding number. For $k_i$ not rational
we get an analytic mess; the direction in the space of B-fields does not
close up to give a circle. 

There is no obvious mirror to these rays since we consider $\Omega$ only
up to scale. So mirror symmetry predicts an isomorphism
between $\K^\C_X$ and the moduli space $\M$ of
complex structures on $\X$, and a distinguished limit in $\M$, the
\emph{large complex structure limit point} (LCLP), the mirror of the
LKLP $(0,\ldots,0)\in\Delta^n$ above. Morrison has given a rigorous definition
of LCLPs and the canonical coordinates on $\M$ dual to the $z_i$ on
$\K_X^\C$; see Section \ref{LCLP} below. The holomorphic curves in
$(\Delta)^n$ described above, corresponding to rational rays of K\"ahler
forms, give degenerations of (the complex structure on) $\X$ to the LCLP
whose monodromy we will discuss in Section \ref{fibs}.

LCLPs play a vital role in mirror symmetry; in fact mirror symmetry is really
a statement about LCLPs and families of Calabi-Yau manifolds near LCLPs. Most predictions
only really hold near or at the LCLP, and complex structure moduli space
only looks like $\Delta^n$
near the LCLP. For instance, manifolds can have many LCLPs and accordingly
many mirrors.
This also explains one obvious paradox -- that \emph{rigid} Calabi-Yau manifolds,
those with no complex structure deformations, $h^{2,1}=0$, and so no LCLP,
can have no mirror,
since a K\"ahler (or symplectic) manifold has $h^2=h^{1,1}\ne0$.

The first predicted refinement of (\ref{hodge}) is, as discussed in
the introduction, that the \emph{variation of Hodge structure} (VHS) on $\X$ should
be describable in terms of \emph{Gromov-Witten} invariants of $X$. Here VHS is
governed by how the ray $\C.\Omega_t=H^{3,0}(\X_t)$ sits inside $H^3(\X_t,\C)$
as the complex structure on $\X_t$ varies, parametrized by $t\in\M$.
By Poincar\'e duality, it is sufficient to know how $\Omega_t$ pairs with
$H_3(\X)$, i.e. to compute the \emph{period integrals}
$$
\int_{A_i}\Omega_t,\qquad i=1,\ldots,2k=2h^{2,1}+2,
$$
where $A_i$ form a basis of $H_3(\X,\Z)$. (In fact we can choose the $A_i$
to be a symplectic basis, $A_i.A_j=\delta_{i+k,j}$, and then knowledge of
only the periods of the first $k$ $A_i$ suffices, locally in moduli space.) These periods determine
$\Omega_t$ and so the Yukawa coupling
\beq{Yuk}
H^1(T\X_t)^{\otimes3}\rt{\cup}H^3(\Lambda^3T\X_t)\Rt{\,\ip\Omega_t^{\otimes2}\!}
H^3(K_{\X_{\!t}})\cong\C.
\eeq
On $X$ we get the cubic form on $H^2(X)$ described roughly in the introduction
in terms of numbers of rational curves in $X$. These numbers are in fact
independent of the almost complex structure on $X$ so long as it is compatible
with the symplectic form $\omega$, so give the symplectic invariants of
Gromov and Witten.
The cubic form depends on $\omega=\omega_t$ as it moves in $\K_{X_t}$ (or
in $\K^\C_{X_t}$, replacing $\omega_t$ by $-i(B_t+i\omega_t)$). Under the predicted
local isomorphism $\K^\C_X\cong\M$ near the LKLP and LCLP, the equality
of these cubic forms gives the predictions of numbers of rational curves
in $X$ mentioned in the introduction. This has been carried out, and the
predictions checked rigorously, in quite some generality, for instance for
mirror pairs produced by Batyrev's toric methods.

There is, of course, a
flat connection, the \emph{Gauss-Manin connection} on the bundle over $\M$
with fiber $H^3(\X_t,\C)$ over $t\in\M$, given by the local system
$H^3(\X_t,\Z)\subset H^3(\X_t,\C)$. Mirror to this, Dubrovin has shown
how to put a flat connection on the bundles with fibers $H^2(X_t)$ and
$H^{\text{ev}}(X_t)$ using Gromov-Witten invariants.

\section{Homological mirror symmetry}

Building on work of Witten, in 1994 Kontsevich \cite{K} proposed a remarkable
conjecture
that purported to explain mirror symmetry, all the more surprising because
it appeared to have little to do with what was thought to be mirror symmetry
at the time. The conjecture is now reasonably well understood,
while the link to Gromov-Witten invariants and Yukawa couplings is more mysterious,
though it is known how both data should be encoded in the conjecture.

Kontsevich proposed that mirror symmetry is a (non-canonical)
equivalence of triangulated categories between the derived Fukaya category $D^{\mathcal F}(X)$ of $(X,\omega)$
and the bounded derived
category of coherent sheaves $D^b(\X)$ on its mirror $\X$. This second category
consists of chain complexes of holomorphic bundles, with quasi-isomorphisms
(maps of chain complexes inducing isomorphisms on cohomology)
formally inverted, i.e. decreed to be isomorphisms. For zero B-field
the first category should be constructed from \emph{Lagrangian submanifolds}
$L\subset X$ carrying flat unitary connections $A$. That is, $L$ is middle
(three) dimensional, and
$$
\omega\res L\equiv 0,\quad F_A=0.
$$
For $B\ne0$ this needs modifying to $F_A+2\pi iB.\id=0$ (so in particular
we require that $L$ satisfies $[B|_L]=0\in H^2(L,\R/\Z)$). There are also
various technical conditions like the choice of a relative spin structure,
the \emph{Maslov class} of $L$ must vanish (i.e. the map $(\Omega|_L/\text{vol}_L)
\colon L\to\C^*$ has winding number zero) and we pick a grading
on $L$ (a choice of logarithm of this map). Morphisms are defined by \emph{Floer
cohomology} $HF^*$ of Lagrangian submanifolds; roughly speaking this assigns a  vector space to each intersection point (the Homs between the fibers of the two unitary bundles carried by the Lagrangians at this point), made into
a chain complex by a certain counting of holomorphic disks between intersection
points. Deep work of Fukaya-Oh-Ohta-Ono shows that this gives the structure
of an $A^\infty$-category which can then be ``derived" into a triangulated
category in a formal way by taking ``twisted cochains". The construction
is still very technical and hard to calculate with, but the key points are
that we
get a category depending only on the symplectic structure, that certain
``unobstructed" Lagrangian submanifolds give objects of this category, and
that \emph{Hamiltonian isotopic} unobstructed Lagrangian submanifolds give
isomorphic objects.

Since the introduction of \emph{D-branes} there is a physical interpretation of
this conjecture in terms of \emph{open string theory};
the objects of the two categories are boundary conditions for open strings,
and morphisms correspond to strings beginning on one object and ending on
the other. So, for instance, intersections of Lagrangians give morphisms
corresponding to constant strings at the intersection point, while the
Floer differential gives instanton tunneling corrections.

One paradox this formulation immediately sheds light on is automorphisms
on both sides of mirror symmetry.
While symplectomorphisms of $(X,\omega)$ are abundant, there
are few holomorphic automorphisms of a Calabi-Yau $\X$. The former induce
autoequivalences of $D^{\mathcal F}(X)$; Kontsevich's suggestion is that mirror
to this there should be an autoequivalence of $D^b(\X)$; this need not be
induced by an automorphism of $\X$. Motivated by this, groups of autoequivalences
of derived categories of sheaves of Calabi-Yau manifolds have now been found that
were predicted by mirror symmetry; a few are mentioned below. Thus homological
mirror symmetry suggests
a SCFT is equivalent to a triangulated category, and the ambiguities in geometrizing
a SCFT (finding a Calabi-Yau of which it is a $\sigma$-model) are seen in the category -- not all automorphisms come from an automorphism
of a Calabi-Yau, Calabi-Yau manifolds $\X$ with equivalent derived categories give multiple mirrors to $X$, and not all appropriate categories need even come
from a Calabi-Yau. Supporting this suggestion, Bondal-Orlov and Bridgeland
have shown that indeed birational Calabi-Yau manifolds $\X$ have equivalent derived
categories.

Finally, Kontsevich explained how deformation theory of the categories should
involve derived morphisms on the product from the diagonal (thought of as a Lagrangian in
the A-model, its structure sheaf as a coherent sheaf in the B-model)
to itself, giving quantum cohomology in the A-model and Hodge structure
in the B-model. For instance the holomorphic disks used to compute the Floer
cohomology of the diagonal on the product $X\times X$ give holomorphic
rational curves on $X$. So one should be able to see some parts of
``classical" mirror symmetry.

Below, as we describe more of the geometry of mirror symmetry that has emerged
since Kontsevich's conjecture, we will mention at each stage how his conjecture
fits in with it.

\section{The SYZ conjecture}

To recover more geometry from Kontsevich's conjecture, there are some
objects of $D^b(\X)$ that obviously reflect the geometry of $\X$ -- the structure
sheaves $\O_p$ of points $p\in\X$. Calculating their self-Homs,
$$
\Ext^*(\O_p,\O_p)\cong\Lambda^*T_p\X\cong\Lambda^*\C^3\cong H^*(T^3,\C),
$$
shows that if they are mirror to Lagrangians $L$ in $X$ (with flat connections
$A$ on them) then we must have
$$
HF^*((L,A),(L,A))\cong H^*(T^3,\C),
$$
as graded vector spaces. Since the left hand side is, modulo instanton corrections,
$H^*(L,\C)^{\oplus r}$, where $r$ is the rank of the bundle carried by $L$,
this suggests that the mirror should be $L\cong T^3$ with a flat $U(1)$ connection
$A$ over it. There are then reasons why the Floer cohomology of such
an object should not be quantum corrected, and so be isomorphic to Ext$^*(\O_p,\O_p)$.

For any Lagrangian $L$, the symplectic form gives an isomorphism between
$T^*L$ and its normal bundle $N_L$; thus Lagrangian tori have trivial
normal bundles, and locally one can fiber $X$ by them. Thus one might
hope that $X$ is fibered by Lagrangian tori, and the mirror $\X$ is (at
least over the locus of smooth tori) the dual fibration. This is because
the set of flat $U(1)$ connections on a torus is naturally the dual torus.

This is the kind of philosophy that led to the SYZ conjecture \cite{SYZ},
though Strominger, Yau and Zaslow were working with physical D-branes, and
not Kontsevich's conjecture. Therefore their D-branes are not the``topological
D-branes" of Kontsevich, but those minimizing some action. That is, instead
of holomorphic bundles in the B-model, we deal with bundles with a compatible
connection satisfying an elliptic PDE (like the \emph{Hermitian-Yang-Mills
equations} (HYM), or some perturbation thereof); instead of Lagrangian
submanifolds up to Hamiltonian isotopy in the A-model, we consider \emph{special
Lagrangians} (sLags) (\ref{slag}). 
The SYZ conjecture is that a Calabi-Yau $X$ should admit a sLag
torus fibration, and that the mirror $\X$ should admit a fibration which
is dual, in some sense.

A sLag is a Lagrangian submanifold of a Calabi-Yau
manifold $X$ satisfying the further equation that
the unit norm complex function (phase)
\beq{slag}
\frac{\Omega|_L}{\vol\_L}=e^{i\theta}=\,\text{constant}\,.
\eeq
(So sLags have Maslov class zero, in particular.) This equation uses the
complex structure on $X$ as well as the symplectic structure, and the resulting
Ricci-flat metric of Yau, to define a metric on $L$ and so its Riemannian
volume form $\vol\_L$. SLags are \emph{calibrated} by $\Re(e^{-i\theta}\Omega)$
and so minimize volume in their homology class. This is similar to the HYM
equations on the mirror
$\X$, which are defined on holomorphic bundles on the complex manifold $\X$
via a K\"ahler form $\omega$, and minimize the Yang-Mills action. The
\emph{Donaldson-Uhlenbeck-Yau theorem} says that for holomorphic
bundles that are \emph{polystable} (defined using $[\omega]$, this is true
for the generic bundle), there is a \emph{unique} compatible
HYM connection. Thus modulo stability,
HYM connections are in one-to-one correspondence with holomorphic bundles.
Thomas and Yau conjecture, and prove in some special cases, a similar correspondence for (special) Lagrangians: that modulo a stability condition
(which can be formulated precisely), sLags are in one-to-one correspondence
with Lagrangian submanifolds up to Hamiltonian isotopy. That is, there should
be a \emph{unique} sLag in the Hamiltonian isotopy class of a Lagrangian if and
only if it is stable. Currently only the uniqueness part of this
conjecture has been worked out, but morally at least, we do not lose much by
considering only \emph{Lagrangian} torus fibrations.

The SYZ conjecture is thought to hold only near the LCLPs and
LKLPs of $X$ and $\X$; away from these the sLag fibers may start to
cross. Due to work of Joyce, the discriminant locus of the fibration on $X$
is expected to be a \emph{codimension one} ribbon graph in a base $S^3$ near
the limit points, while the discriminant locus of the dual fibration
$\X$ may be different -- i.e. the smooth parts of the fibration
and its dual are compactified in different ways.
In the limit of moving to the limit points, however, both discriminant loci
shrink onto the same codimension two graph. In this limit the fibers shrink to zero size, so that $X$ (with
its Ricci-flat metric) tends in the Gromov-Hausdorff sense to its base
$S^3$ (with a singular metric). This formal picture has been made precise
in two dimensions, for $K3$-surfaces, by Gross and Wilson. The limiting
picture suggests that if we are only interested in topological
or \emph{Lagrangian} torus fibrations then we might hope for codimension two
discriminant loci, and such fibrations might make sense well away from
limit points. Work of Gross and Ruan carries this out in examples such
as the quintic and its mirror, and makes sense
of dualizing the fibration by dualizing monodromy around the discriminant
locus and specifying a canonical compactification over the discriminant
locus. This gives the correct topology for toric varieties and their mirrors,
and flips Hodge numbers (\ref{hodge}), for instance. Approaching the LCLP
in a different way (in the example (\ref{ni}) this corresponds to altering
the rational numbers $k_i$) can give a different graph and different fibration
on $X$; the dual fibration can then be a topologically different manifold,
giving a different birational model of the mirror $\X$.

We focus only on Lagrangian fibrations, as they are better behaved and understood.
We can expect them to be $C^\infty$ fibrations with codimension two discriminant
loci, for instance. Below we see how to put a complex structure on the
smooth part
of the fibration, but extending this over the compactification is much
harder and will involve ``instanton corrections" coming from holomorphic
disks. Fukaya \cite{F} has beautiful conjectures about this that will
explain a great deal more of mirror symmetry, but they will not be discussed
here.

\section{Lagrangian torus fibrations} \label{fibs}

If $(X^{2n},\omega)\rt{\pi}B^n$ is a \emph{smooth} Lagrangian fibration with
compact
fibers, then the fibration is naturally an affine bundle of torus groups
(i.e. a bundle of groups once we pick a Lagrangian zero-section -- an identity
in each fiber), and the base $B$ inherits a natural integral affine structure:
it looks like a vector space $V$ with an integral structure
$V\cong\Lambda\otimes_\Z\R$ up to translation by elements of $V$. This
is the classical theory of \emph{action-angle variables}. $T^*_bB$
acts on the fiber $X_b=\pi^{-1}(b)$: by pullback and contraction with the
symplectic form, $\sigma\in T^*_bB$ gives a vector field $\underline{\sigma\!}\,$
tangent to $X_b$, and the time-one flow along $\underline{\sigma\!}\,$
gives the action. By compactness
and smoothness of $X_b$ the kernel is a full rank lattice $\Lambda_b\subset
T^*_bB$, giving the isomorphism
$$
X_b\cong T^*_bB/\Lambda_b.
$$
We define the integral affine structure on $B$ by specifying the integral
affine functions $f$ (up to translation) to be those whose time-one flow
along $\underline{df\!}\,$ is the identity (i.e. on the universal cover the
time-one flow is to a section of the bundle of lattices $\Lambda$).

The situation that interests us is where $B$ is a 3-manifold
$\overline B$ (usually $S^3$) minus a graph; then the monodromy around
the graph preserves the integral affine structure:
\beq{intaff}
\pi_1(B)\to\R^3\rtimes GL(3,\Z).
\eeq
A great deal of mirror symmetry can be seen from just this knowledge of
the smooth locus of the fibration; in particular Gross \cite{G} has
shown how mild assumptions about the compactification (with singular fibers
over $\overline{B}\take B$) are enough to determine much of the topology
of $X$.
The dual fibration $\ppi$ should have the monodromy \emph{dual}
to (\ref{intaff}), and he shows how this implies
the switching of the Hodge numbers (\ref{hodge}) by the Leray spectral
sequence; the rough idea being the obvious isomorphism
$$
R^i\pi_*\R\cong\Lambda^iTB\cong\Lambda^{3-i}T^*B\cong R^{3-i}\ppi_*\R
$$
induced by a trivialization of $\Lambda^3TB$.
That is, intuitively, the flipping of Betti numbers arises by representing cycles by those with linear intersection
with the fibers, and replacing this linear space by its annihilator in the
dual torus. This also agrees with the equivalence taking Lagrangians to coherent
sheaves described in the next section.

The dual fibration $\ppi$ has a natural complex structure; here the
affine structure is essential, as in general a tangent bundle $TB$ only
has a natural almost complex structure along its zero-section. Since, up
to translation, locally $B\cong V$ is a vector space, $TB\cong V\times V\cong
V\otimes_\R\C$ has a natural complex structure which descends to
\beq{dualfib}
\ppi\colon\X=TB/\Lambda^*\to B.
\eeq
Gross suggests that the \emph{B-field} on $X$ should lie in the piece
$$
H^1(R^1\pi_*\R/\Z)=H^1(TB/\Lambda^*)
$$
of the Leray spectral sequence for $H^2(X,\R/\Z)$. So it gives a \v{C}ech cocycle $e$ on overlaps of an open cover of
$B$ with values in the dual bundle of groups $TB/\Lambda^*$. Using this
to twist (\ref{dualfib}) and reglue it via transition functions translated
by $e$, we get a new complex manifold ($e$ is locally constant so
translation by $e$ is holomorphic)
which we think of as mirror to $X$ with complexified form $B+i\omega$.
In this way Gross manages to match up complexified symplectic deformations
of $X$ with complex structures on $\X$.

\section{The two-torus}

Mirror symmetry is non-trivial even for the simplest Calabi-Yau -- the two-torus.
We write it as an SYZ fibration $T^2\rt{\pi}B=S^1$, and write $B$ as
$\R/a\Z$ with its standard integral affine structure induced by $\Z\subset\R$.
This trivializes $T^*B=B\times\R$ and the lattice $\Lambda$ in it as $B\times\Z
\subset B\times\R$. So as a symplectic manifold,
\beq{rectangle}
T^2=\frac{T^*S^1}\Lambda=\frac{[0,a]\times[0,1]}{\text{\small$(0,p)\!\sim\!
(a,p),\ (q,0)\!\sim\!(q,1)$}}\,,
\eeq
with symplectic coordinates $(q,p)$ in which the symplectic form is $\omega=
dp\wedge dq$ (so $\int_{T^2}\omega=a$). Again the B-field, $b\in H^1(R^1\pi_*
\R/\Z)=H^2(T^2,\R/\Z)$, is in $H^1$ of the locally constant
sections of the dual fibration.

In our trivialization $B\cong\R/a\Z$, $\Lambda^*\subset TB$ is also standard:
$B\times\Z\subset B\times\R$, so the mirror has the same description
(\ref{rectangle}) in which the complex structure is standard: $J\partial_p=
\partial_q$. That is, $p+iq$ gives a local holomorphic coordinate.

For nonzero B-field $b\ne0$, twisting the dual fibration by $b$ gives
\beq{oblong}
T^2=\frac{T^*S^1}\Lambda=\frac{[0,a]\times[0,1]}{\text{\small$(0,p)\!\sim\!
(a,b+p),\ (q,0)\!\sim\!(q,1)$}}\,,
\eeq
again with holomorphic structure given by $p+iq$ and SYZ fibration $\ppi$
being projection onto $q$. So as a complex manifold the mirror is $\C$ divided
by the lattice
$$
\Lambda=\langle 1,b+ia\rangle.
$$
Changing $b$ to $b+1$ does not alter this lattice, so the construction
is well defined for $b\in\R/\Z\cong H^1(R^1\pi_*\R/\Z)$, and we have the
standard description of an elliptic curve via its period point $\tau=b+ia$
in the upper half plane (as $a>0$). Mirror symmetry has indeed swapped
the complexified symplectic parameter $b+ia=\int_{T^2}(b+i\omega)$ for
the complex structure modulus $\tau=b+ia$. $SL(2,\Z)$ acts on both sides
(in the standard way on $\tau$, and as symplectomorphisms modulo those
isotopic to the identity on the A-side) permuting
the choices of SYZ fibration. We note that in this case, the fibrations
are \emph{special} Lagrangian in the flat metric, with no singular fibers.

Polishchuk and Zaslow have worked out in detail how Kontsevich's conjecture
works in this case. The general picture for any torus fibration is an extension
of the fiberwise duality that led to SYZ. Namely, Lagrangian multisections
$L$ of the fibration, of degree $r$ over the base, give $r$ points on each
fiber, and so $r$ flat $U(1)$ connections on the dual fiber. The resulting
$U(1)^{\times r}$ connections can be glued together and
twisted by the flat connection on $L$, to give a rank $r$ vector bundle
with connection on the mirror. Arinkin and Polishchuk show that in general
the Lagrangian condition implies the integrability condition $F^{0,2}=0$
of the resulting connection, giving a holomorphic structure on the bundle.
Leung-Yau-Zaslow show that the special Lagrangian condition gives a perturbation
of the HYM equations on the connection. Branching of sections has been
dealt with by Fukaya, and requires instanton corrections from holomorphic
disks. Other Lagrangians with linear intersection with the fibers can be
dealt with similarly. $T^2$ is simpler because all Lagrangians with vanishing
Maslov class can be isotoped into straight lines (i.e. sLags
in the flat metric) with no branching. The upshot is that the slope of
the sLag over the base corresponds to the slope $\left(\int_{T^2}c_1/
\rk\right)\in[-\infty,\infty]$ of the mirror sheaf.

\subsection*{The large complex structure limit}

The LKLP for $T^2$ is clearly $\lim a\to\infty$. On the mirror then, the
LCLP is at $\tau=b+ia\to b+i\infty$, the nodal torus compactifying the
moduli of elliptic curves.
\emph{Metrically}, however, in the (Ricci) flat metric, things look
different; if we rescale to have fixed diameter, the torus collapses to
the base of its SYZ fibration, and \emph{all} of its fibers contract. This is
an important general feature of the difference between complex and metric
descriptions of LCLPs; see the description of the quintic in the next section.

We note that as in the compactifications of Section \ref{1}, the monodromy
around this LCLP is given by rotating the B-field: $b\mapsto b+1$. This
gives back the same elliptic curve, but after a monodromy diffeomorphism
$T$ which we see from (\ref{oblong}) to be
$$
T\colon\ q\mapsto q, \quad p\mapsto p+q/a.
$$
On $H^1(T^2)=\Z[\text{fiber}]\oplus\Z[\text{section}]$
this acts as
\beq{mat}
T_*=\mat{1}{1}{0}{1}.
\eeq
This is called a \emph{Dehn twist}. Picking the zero-section $O=\{p=0\}$
in the mirror (\ref{oblong}) when $b=0$, this is taken to the section
$$
T(O)=\{p=q/a\},
$$
and $T$ is in fact translation by this section $T(O)$ on $T^2$, using the
group structure on the fibers (now we have chosen a zero-section). Again,
Gross \cite{G} has shown that this is a general feature of LCLPs.

If we pick a K\"ahler structure on this family of complex tori, $T$ turns
out to
be a symplectomorphism. Importantly, its mirror is \emph{not} a holomorphic
automorphism, but an equivalence of the derived category of coherent sheaves.
As above, the section $T(O)$ corresponds to a slope-one line bundle $L$
on the mirror, and the monodromy action corresponds to
\beq{L}
\otimes L\colon D^b\to D^b
\eeq
on the
derived category. Again, this is a more general feature of these LCLPs,
with $L$ such that $c_1(L)$ equals the symplectic form which generated the
ray along which the original LKLP was reached. In general the SYZ fiber is
the invariant cycle under $T_*$ (\ref{mat}), and, on the mirror, structure
sheaves of points are invariant under $\otimes L$. On the cohomology of
$T^2$, cupping with $ch(L)=e^{c_1(L)}=1+c_1(L)$ has the same action (\ref{mat})
on $H^{\text{ev}}=\Z(c_1(L))\oplus\Z(1)$.

Notice we have used the choices of fibration and zero-section to produce
the equivalence of triangulated categories and to equate the monodromy
actions. Kontsevich's conjectural equivalence is not canonical, but is
fixed by a choice of fibration and zero-section. In turn a fibration should
be fixed by a choice of LCLP or LKLP from the resulting collapse (in the
Ricci flat metric) onto a half dimensional $S^n$ base. The choice of zero-section is then rather arbitrary (as monodromy about the LCLP changes it)
but determines the equivalence of categories. Different choices of section
give different equivalences, differing for instance by the monodromy transformation
$\otimes L$ (\ref{L}).

Another point of view is that a Lagrangian fibration and zero-section determine a
group structure on the fibers and so on the Fukaya category (translating
Lagrangian multisections by multiplication on each fiber). This corresponds
to a choice of tensor product on the derived category of the mirror; the
identity for this product is then the structure sheaf $\O_X$ mirror to
the zero-section, and an ample line bundle is given by the action of the
monodromy transformation $L=T(\O_X)$; $T$ then acts as $\otimes T(\O_X)$
(\ref{L}). Since $X$ is determined by the graded ring
$$
\bigoplus_{j\gg0}H^0_X(L^j)=\bigoplus_{j\gg0}\Hom^*(\O_X,T^j(\O_X)),
$$
one might also try to construct $X$ purely from the zero-section $O$ and
LCLP monodromy on $\X$, as
$$
X=\Proj\bigoplus_{j\gg0}HF^*(O,T^j(O)).
$$
A problem is to show that $\oplus_{j\gg0}HF^0(O,T^j(O))$ is finitely generated;
a related problem is to show that for $j\gg0$, the above Floer
homologies vanish except for $*=0$.

We now turn to the quintic, where we will see how to identify the (homology
classes of the) zero-section and fiber in general using Hodge theory.

\section{The quintic 3-fold}

The simplest Calabi-Yau 3-fold is given by the zeros $Q$ of a homogeneous quintic polynomial
on $\PP^4$, i.e. an anticanonical divisor of $\PP^4$. By adjunction this has
trivial canonical bundle, and so is Calabi-Yau. By the Lefschetz hyperplane
theorem, it has $h^{1,1}=1$,
so computing its Euler number to be $e=-200$, we find $h^{2,1}=101$ is its
number of complex deformations. Alternatively this can be seen by showing
all such deformations are themselves quintics, then dividing the 126-dimensional
space of quintic polynomials by the 25-dimensional $GL(5,\C)$. Thus its
mirror has 1 complex structure deformation and 101 K\"ahler classes.

Greene and Plesser prescribed the following mirror. Take the special one-dimensional
family of Fermat quintics
\beq{GP}
Q_\lambda=\left\{\sum_{i=0}^4x_i^5-\lambda\prod_{i=0}^4x_i=0\right\}\subset\PP^4,
\eeq
with the action of $\{(\alpha_0,\ldots,\alpha_4)\in(\Z/5)^5\colon
\prod_i\alpha_i=1\}\cong(\Z/5)^4$ given by rescaling the $x_i$ by fifth roots
of unity. Dividing by the diagonal $\Z/5$ projective stabilizer we get a
free $(\Z/5)^3$-action; the mirror
of the quintic is any crepant ($K=\O$) resolution of the quotient:
$$
\check Q_\lambda=\widehat{\frac{Q_\lambda}{(\Z/5)^3}}.
$$
Different resolutions have isomorphic $H^2$ but different K\"ahler cones therein. The union of these K\"ahler cones is
the moveable cone, whose complexification is locally isomorphic
to the complex structure moduli space of $Q$. $h^{1,1}(\check Q_\lambda)=101$
for any crepant resolution, and $h^{2,1}(\check Q_\lambda)=1$ corresponds
locally
to the one complex structure deformation (\ref{GP}). In fact for $\alpha^5=1$,
multiplying $x_0$ by $\alpha$ shows that $\check Q_\lambda\cong\check
Q_{\alpha\lambda}$, and $\lambda^5$ parametrizes the complex structure moduli.

The large complex structure limit point is at $\lambda=\infty$, i.e. it
is a resolution of the quotient of the union of hyperplanes
\beq{lclp}
Q_\infty=\left\{\prod_{i=0}^4x_i=0\right\}=\{x_0=0\}\cup\ldots\cup\{x_4=0\}.
\eeq
This is a union of toric varieties, each with a $T^3$-action inherited from the toric $T^4$-action
on $\PP^4$. Much more generally, Batyrev's construction considers the
anticanonical divisors (and even more generally, complete intersections)
in toric varieties fibered over the boundary of the
moment polytope, and takes as mirror the anticanonical divisor of the toric
variety associated to the \emph{dual} polytope. Most of the geometry
is visible in this quintic example, however.

(\ref{lclp}) is the analogue of the nodal torus of the last section, and
we emphasize
again that metrically it looks nothing like this; the Ricci-flat
metric collapses the $T^3$ toric fibers to the base $S^3$ (with a singular
metric). General LCLPs look rather similar, with such ``as bad as possible"
normal crossing singularities. Smoothing a local model (in $x_0=1$)
$\prod_{i=1}^4 x_i=0$ we can see the tori in $\big\{\prod_{i=1}^4 x_i
=\epsilon\big\}$:
\beq{T3}
T^3\!=\!\left\{|x_1|=\delta_1,|x_2|=\delta_2,|x_3|=\delta_3,
x_4=\frac\epsilon{x_1x_2x_3}\right\}\!.
\eeq
These are even Lagrangian in the standard symplectic form on the local model,
and fiber the smoothing over the base $\{(\delta_1,\delta_2,\delta_3)\}$.
It turns out that metrically, these tori which vanish into the normal
crossings singularity at the LCLP actually form a large part of the smooth
Calabi-Yau. This enlightens the apparent paradox between the SYZ conjecture and the
Batyrev construction, i.e. why a vertex of the original moment polytope
(corresponding to the deepest type of singularity $(0,0,0,0)\in
\big\{\prod_{i=1}^4x_i=0\big\}$)
can be replaced by the dual 3-dimensional face in the dual polytope. This
was first suggested by Leung and Vafa.

Gross and Siebert \cite{GS} exploit this to extend SYZ and Batyrev's
construction to non-toric LCLP Calabi-Yau manifolds;
it is only the local toric nature of the normal crossing singularities of
the LCLP that they use. It seems possible that their construction will give
the mirrors of all Calabi-Yau manifolds with LCLPs. Much of mirror symmetry should
soon
be reduced to graphs (the discriminant locus of a Lagrangian torus fibration)
in spheres, and further graphs over which D-branes (such as holomorphic
curves) fiber, as in recent conjectures of Kontsevich and Soibelman and Fukaya
\cite{F}. Before long
it may be possible to write down a triangulated category in terms of
such data. The full geometric story (involving Joyce's description of sLag
fibrations, for instance) is still some way off, however; we cannot
even write down an explicit Ricci-flat metric on a compact Calabi-Yau.

\subsection{Monodromy around the LCLP} \label{LCLP}

As well as the SYZ torus fiber (\ref{T3}) we can also see a Lagrangian
zero-section
on the quintic and its mirror as a component of the real locus of (\ref{GP})
for $\lambda>5$.
Remarkably, like the torus (\ref{T3}), this cycle was already described and
used in \cite{COGP}, long before the relevance of torus fibrations was suspected.

Gross and Ruan are able to describe the quintic and its mirror (at least
topologically or symplectically) very explicitly as a simple torus fibration
over this $S^3$ with a natural integral affine structure and codimension
two graph discriminant locus. See for example \cite{GHJ}.

Under monodromy about $\lambda=\infty$, the zero-section is moved to another
section $T(O)$, and $T$ is given by translation by $T(O)$ using the group
structure on the fibers. This is the analogue of the Dehn twist (\ref{mat}),
and one can choose a basis of $H_3(\check Q)$ (with first element the invariant
cycle, the $T^3$-fiber, second element a cycle fibered over a curve in $S^3$,
third fibered over a surface, and last the zero-section itself) such that
\beq{matrix}
T_*=\left(\begin{array}{cccc} 1 & 1 & * & * \\ 0 & 1 & * & * \\ 0 & 0 &
1 & * \\ 0 & 0 & 0 & 1
\end{array}\right).
\eeq
Like the Dehn twist (\ref{mat}), it turns out that $T_*$ is \emph{maximally
unipotent}; that is, in $n$-dimensions, we have
$$
(T_*-1)^{n+1}=0 \quad\text{but}\quad (T_*-1)^n\ne0.
$$
Again, this is a general feature of LCLPs as formulated by Morrison \cite{M}
as part of the definition.

This should be compared with the Lefschetz operator $L=\cup\,\omega$ on the
cohomology of the mirror, which also satisfies $L^n\ne0,\ L^{n+1}=0$
(or, more relevantly, $\exp(L)$, which satisfies $(e^L-1)^n\ne0,\ 
(e^L-1)^{n+1}=0$). Their
similarity was noticed by the Griffiths school working on VHS in the late
1960s! Now we know that for Calabi-Yau manifolds at a LCLP dual
to a LKLP along a ray $\omega=c_1(L)$ on the mirror, they should
be considered mirror operators (up to some factors of the Todd class of
the underlying Calabi-Yau, to do with the relationship between the Chern
character $e^\omega$ of the line bundle $L$ (\ref{L}) and the Riemann-Roch
formula).

Both, by linear algebra of the nilpotent operator $N=\log T_*=\sum_{k=1}^n(T_*-1)^k$,
induce a natural filtration $W_\bullet\colon 0\le W_0\le\ldots\le W_{2n}=H$
on the cohomology on which they operate (which is $H=H^n$ for $N=\log T_*$
and $H=H^{\text{ev}}$ for $N=L=\cup\,\omega$):
\beq{filt}
\!\!\!\!\!\! 0\le\im(N^n)\le\im(N^{n-1})\cap\ker(N)\le\ldots \qquad
\eeq $$
\quad\qquad \le\ker(N^{n-1})+\im(N)\le\ker(N^n)\le H.
$$
We refer to the references for the construction of this \emph{monodromy weight
filtration}. It plays a key role in studying
degenerations of varieties and Hodge structures, in this case as we approach
the LCLP. It is a beautiful result of Gross that this filtration coincides
with the Leray filtration on $H^n$ induced by the fibration. That is, under
Poincar\'e duality, the weight
filtration on cycles is by the minimal dimension (over all homologous
cycles) of the image in the base over which the cycle is fibered.
So the first graded piece is spanned by the invariant cycle, the $T^3$-fiber,
supported over a point, and the last by the zero-section; cf. (\ref{matrix}).
(Similarly on the mirror, the filtration for the Lefschetz operator $\cup\,
e^\omega$ has first piece spanned by the cohomology class of a point, which
is invariant under the monodromy action $\otimes L$ (\ref{L}), etc.)

Letting $\gamma_0$ be the class of a fiber and $\gamma_1$ span
$W_2/W_0$ (which is one dimensional) over the integers, then
$T_*\gamma_1=\gamma_1+\gamma_0$. It follows that
$$
q=\exp\left(2\pi i\,\frac{\int_{\gamma_1}\Omega}{\int_{\gamma_0}\Omega}\right)
$$
is invariant under monodromy. This is the higher dimensional analogue of
the coordinate $\exp(2\pi i\tau)$ on the moduli space of elliptic curves,
where $\tau$ is the period point. It is this coordinate $q$ that is mirror
to the coordinate
$$
\int_{\text{line}}\omega
$$
on the K\"ahler moduli space on the mirror quintic, which allows one to
compute the correspondence between VHS and Gromov-Witten
invariants mentioned in the introduction.

More generally, following Morrison \cite{M},
one can make a rigorous definition of a LCLP using features
we have seen above, extended to the case of $h^{2,1}>0$; see for instance
\cite{CK}. Roughly, the upshot is that $\M$ (of dimension $s=h^{2,1}(\X)$)
should be compactified with $s$ divisors $(D_i)_{i=1}^s$ (parameterizing
singular varieties) forming a normal crossings divisor meeting at the LCLP,
with monodromies $T_i$ about them. There should be
a unique (up to multiples) integral cycle $\gamma_0$ (our torus fiber, of
course) invariant under all $T_i$, and cycles $(\gamma_i)_{i=1}^s$ such
that
$$
\tau_i=\frac{\int_{\gamma_i}\Omega}{\int_{\gamma_0}\Omega}
$$
is logarithmic at $D_i$; i.e. $\tau_i=\frac1{2\pi i}\log(z_i)$, where
$z_i$ is a local parameter for $D_i=\{z_i=0\}$.

So $z_i=\exp(2\pi i\tau_i)$ form local coordinates for moduli space, mirror
to the polydisk coordinates (\ref{q})
on $\K^\C_X$. The direction we approached the LKLP in that section corresponds
to the holomorphic curve $z_i^{k_j}=z_j^{k_i}$ (\ref{ni}) we take through
the LCLP $(z_i=0\ \forall i)$, and the monodromy $\sum N_iT_i$ varies accordingly,
but the corresponding weight filtration $W_\bullet$ remains constant if
$k_i\ne0\ \forall i$, by a theorem of Cattani and Kaplan.

Morrison then requires that the $(\gamma_i)_{i=0}^s$
should form an integral basis for $W_2=W_3$ (with $\gamma_0$ a basis of $W_0=W_1$).
Finally, part definition and part conjecture, we should be able to choose
that they satisfy $\log T_i(\gamma_j)=\delta_{ij}\gamma_0$.

Of course, as has been emphasized, Morrison's definition of a LCLP is really
where the mathematics and geometry of mirror symmetry begin,
and should have been the starting point of this article. But that would have
required a lot of knowledge of abstract VHS that
are best understood, in this context, through the new geometry of Lagrangian
torus fibrations that mirror symmetry has inspired.

\end{document}